\documentclass[11pt,a4paper,reqno]{amsart} 

\usepackage{amssymb,graphicx}

\usepackage{amssymb}

\newcommand{\koniec}{\begin{flushright}  $\Box $ \end{flushright}}
\newtheorem{theo}{Theorem}[section] 
\newtheorem{prop}[theo]{Proposition}

\newcounter{mnotecount}[section]

\renewcommand{\themnotecount}{\thesection.\arabic{mnotecount}}

\newcommand{\mnote}[1]
{\protect{\stepcounter{mnotecount}}$^{\mbox{\footnotesize
$
\bullet$\themnotecount}}$ \marginpar{
\raggedright\tiny\em
$\!\!\!\!\!\!\,\bullet$\themnotecount: #1} }

\newcommand{\R}{\mathbb{R}}

\newcommand{\Rho}{\mathrm{P}}
\def\p{\partial}
\def\be{\begin{equation}}
\def\ee{\end{equation}}

\def\bea{\begin{eqnarray}}
\def\eea{\end{eqnarray}}

\begin{document}
\date{May 27, 2016}
\title[First integrals and Hamiltonian systems]{First integrals of affine connections and Hamiltonian systems of hydrodynamic type}
\author{}
\author{Felipe Contatto }
\author{Maciej Dunajski}
\address{Department of Applied Mathematics and Theoretical Physics\\ 
University of Cambridge\\ Wilberforce Road, Cambridge CB3 0WA\\ UK.}
\email{felipe.contatto@damtp.cam.ac.uk, m.dunajski@damtp.cam.ac.uk}
\begin{abstract} 
We find necessary and sufficient conditions for a local geodesic flow of an
affine connection on a surface to admit a linear first integral. The conditions are expressed in terms of two scalar 
invariants of differential orders 3 and 4 in the connection. We use this result to find explicit obstructions to the existence of a Hamiltonian formulation of Dubrovin--Novikov type for a given one--dimensional system of hydrodynamic type. We give 
several examples including Zoll connections, and  Hamiltonian systems arising 
from two--dimensional Frobenius manifolds.
 
\end{abstract}   
\maketitle
\section{Introduction}
The existence of a first integral of a geodesic flow of an affine connection puts restrictions on the form of the connection. A generic
connection admits no first integrals. If the connection arises from a
metric, and the first integral is linear in velocities, then the metric admits a one--parameter group of isometries generated by a Killing vector field. Characterising metrics which admit Killing vectors by local tensor obstructions is a classical problem which goes back at least to Darboux \cite{darboux}, and can be solved completely in two dimensions. The analogous characterisation of 
non--metric affine connections has not been carried over in full\footnote{The remarkable exception is the paper of Levine \cite{levine} and its extension \cite{thompson} where the necessary condition for the existence of a first integral was found, albeit
not in a form involving the Schouten and Cotton tensors. The sufficient conditions found in \cite{levine} are not all independent. Levine gives seven tensor conditions, where in fact two scalar conditions suffice.}. It is given in
Theorem \ref{theo_1}, where we construct two invariant scalar obstructions
to the existence of a linear first integral. A non--metric connection can (unlike a Levi--Civita connection) admit precisely two independent linear local first integrals. This case will also be characterised
by a tensor obstruction.
\vskip5pt
As an application of our results we shall, in \S \ref{sec_hydro}, characterise one--dimensional systems of hydrodynamic type which admit a Hamiltonian formulation of the Dubrovin--Novikov type 
\cite{dubrovin}. The existence of such formulation leads to an over--determined system of PDEs, and  we shall show (Theorem \ref{theo_2}) that
this system is equivalent to a condition that  a certain non--metric
affine connection admits a linear first integral. This, together with
Theorem \ref{theo_1} will lead to a characterisation of Hamiltonian, bi--Hamiltonian and tri--Hamiltonian systems of hydrodynamic type.
In \S \ref{sec_ex} we shall give examples of connections resulting from 
hydrodynamic type systems. In particular we shall show that systems arising 
           from
two--dimensional Frobenius manifolds are tri--Hamiltonian.

\vskip5pt
In the remaining part of the Introduction we shall state our main results.
Let $\nabla$ be a torsion-free  affine connection of differentiability class
$C^4$ on a simply connected orientable surface $\Sigma$ (so we require the transition functions of $\Sigma$ to be of class at least $C^6$). A curve
$\gamma:\R\rightarrow \Sigma$ is an affinely parametrised geodesic if
$\nabla_{\dot{\gamma}}\dot{\gamma}=0$, or equivalently
if
\be
\label{flow}
\ddot{X}^a+\Gamma^{a}_{bc}\dot{X}^b\dot{X}^c=0, \quad a, b, c =1, 2
\ee
where $X^{a}=X^{a}(\tau)$ is the curve $\gamma$ expressed in local coordinates 
$X^a$ on an open set $U\subset \Sigma$, $\tau$ is an affine parameter, $\Gamma_{ab}^c$ are the Christoffel symbols of $\nabla$, and we use the summation convention. 
A linear function on $T\Sigma$ given by $\kappa=K_a(X)\dot{X}^a$ is called a 
{\em first integral} if $d\kappa/d\tau=0$ when (\ref{flow}) holds, or 
equivalently if
\be
\label{killing}
\nabla_{(a} K_{b)}=0.
\ee
The following Theorem  gives local necessary and sufficient conditions for a connection
to admit one, two or three linearly independent solutions to the Killing 
equation (\ref{killing}). The necessary conditions involve vanishing of 
obstructions $I_N$ and $T$ given by (\ref{i_n}) and (\ref{tensor_obs}) - 
for these to make sense the connection needs to be at least three times differentiable. 
\begin{theo}
\label{theo_1}
The necessary condition for a $C^4$ 
torsion--free affine connection 
$\nabla$ on  a surface $\Sigma$ to admit a linear first
integral is the vanishing, on $\Sigma$, of 
invariants $I_N$ and $I_S$ given by (\ref{i_n}) and (\ref{i_s}) respectively. For any point $p\in \Sigma$ there exists
a neighbourhood $U\subset \Sigma$ of $p$ such that conditions $I_N=I_S=0$ on $U$
are sufficient for the existence of a first integral on $U$.
There exist precisely two independent linear first integrals on $U$ 
if and only if the tensor  $T$ given by (\ref{tensor_obs}) vanishes and the skew part of the Ricci tensor of $\nabla$ is non--zero on $U$. 
There exist three independent first integrals on $U$  if and only if the connection is projectively flat and its Ricci tensor is symmetric.
\end{theo}
This Theorem will be established by constructing (Proposition \ref{tractor_prop})
a prolongation connection $D$ on a rank--3 vector bundle $\Lambda^1(\Sigma)\oplus\Lambda^2(\Sigma)$ for the over--determined system (\ref{killing}), and restricting the holonomy of its curvature when one, two or three parallel sections exist.  In Proposition \ref{prop_rank1} we shall find all local normal forms of connections from Theorem \ref{theo_1} which admit precisely two linear first integrals. 

Finally we shall consider one--dimensional systems of hydrodynamic type.
Any such system with two dependent variables $(X^1, X^2)$ and two independent variables $(x, t)$ can be written in the so--called Riemann invariants as
\be
\label{hydro_system}
\frac{\p X^1}{\p t}=\lambda^{1}(X^1, X^2)\frac{\p X^1}{\p x}, \quad 
\frac{\p X^2}{\p t}=\lambda^{2}(X^1, X^2)\frac{\p X^2}{\p x},
\ee
where $\lambda^1\neq \lambda^2$ at a generic point.
This system admits a Hamiltonian formulation of the Dubrovin--Novikov type, if it can be written as
\be
\label{hami}
\frac{\p X^a}{\p t}=\Omega^{ab}\frac{\delta H}{\delta X^b},
\ee
where $H[X^1, X^2]=\int_{\R} {\mathcal H}(X^1, X^2) dx$ is the Hamiltonian of hydrodynamic type, and 
the Poisson structure on the space of maps $\mbox{Map}(\R, U)$ is given by
\[
\Omega^{ab}=g^{ab}\frac{\p}{\p x}+{b^{ab}_c}\frac{\p X^c}{\p x}.
\]
The Jacobi identity imposes severe constraints on $g(X^a)$ and $b(X^a)$ -- see
Section \ref{sec_hydro} for details.
We shall prove
\begin{theo}
\label{theo_2}
The hydrodynamic type system (\ref{hydro_system}) admits one, two or three Hamiltonian formulations
with  hydrodynamic Hamiltonians if and only if 
the affine torsion--free connection $\nabla$ defined by its
non--zero components
\begin{eqnarray}
\label{connection_{AB}}
&&\Gamma_{11}^1=\p_1 \ln{A}-2B, \quad \Gamma_{22}^2=\p_2\ln{B}-2A, \quad \Gamma_{12}^1=-\Big(\frac{1}{2}\p_2 \ln{A}+A\Big), \quad
\Gamma_{12}^2=-\Big(\frac{1}{2}\p_1 \ln{B}+B\Big),\nonumber\\
&&\mbox{where}\quad
A=\frac{\p_2 \lambda^1}{\lambda^2-\lambda^1}, \quad
B=\frac{\p_1 \lambda^2}{\lambda^1-\lambda^2}, \quad\mbox{and}\quad 
\p_a=\p/\p X^a
\end{eqnarray}
admits one, two or three  independent linear first integrals respectively.
\end{theo}
This Theorem, together with Theorem \ref{theo_1} leads to explicit
obstructions for the existence of a Hamiltonian formulation (\ref{hami}).
\subsubsection*{Acknowledgements}
The problem of characterising Hamiltonian systems of hydrodynamic type
was suggested to the second author by Jenya Ferapontov in 2008.
F.C. is grateful for the support of Cambridge Commonwealth, European $\&$ International Trust and CAPES Foundation Grant Proc. BEX 13656/13-9. We thank 
Robert Bryant for discussions about invariants of ODEs, and 
Guido Carlet for discussions
about Frobenius manifolds, and for bringing reference \cite{romano} to our attention.  
\section{Killing operator for affine connection}
Given an affine connection $\nabla$ on a surface $\Sigma$, 
its curvature  is defined by
\[
[\nabla_a, \nabla_b]X^c={{R_{ab}}^c}{}_d X^d.
\]
In two dimensions the projective Weyl tensor vanishes,
and the curvature  can be  uniquely decomposed as
\be
\label{formula_for_rho}
{{R_{ab}}^c}{}_d= 
\delta_a{}^c\Rho_{bd}-\delta_b{}^c\Rho_{ad} +B_{ab}\delta_d{}^c,
\ee
where $\Rho_{ab}$ is the Schouten tensor related to the Ricci tensor $R_{ab}={{R_{ca}}^c}_b$
of $\nabla$ by
$\Rho_{ab}=(2/3)R_{ab}+(1/3)R_{ba}$, and 
$B_{ab}=\Rho_{ba}-\Rho_{ab}=-2\Rho_{[ab]}$.
We shall assume that $\Sigma$ is orientable, and choose a volume form $\epsilon_{ab}$ on $\Sigma$. 
We shall also introduce $\epsilon^{ab}$ such that
$\epsilon^{ab}\epsilon_{cb}=\delta_c^a$. These skew-symmetric 
tensors are used to raise and lower indices according to
$V^a=\epsilon^{ab}V_b$ and $V_a=\epsilon_{ba}V^b$.
Then
\[
\nabla_a\epsilon_{bc}=\theta_a\epsilon_{bc},
\]
where $\theta_a=(1/2)\epsilon^{bc}\nabla_a\epsilon_{bc}$. 
Set
$\beta=B_{ab}\epsilon^{ab}$.
\begin{prop}
\label{tractor_prop}
There is a one--to--one correspondence between solutions to the
Killing equations (\ref{killing}), and parallel sections
of the prolongation connection $D$ on a rank--three vector 
bundle $E=\Lambda^1(\Sigma)\oplus\Lambda^2(\Sigma)\rightarrow \Sigma$ defined by
\be
\label{tractor_con}
{\quad{{D}}_a \left(\begin{array}{c}
K_b\\ 
\mu
\end{array} \right)= 
\left(\begin{array}{c} \nabla_a K_b-\epsilon_{ab}\mu \\ 
\nabla_a\mu -\Big({\Rho^b}_a+\frac{1}{2}\beta{\delta^b}_a\Big)K_b+\mu\theta_a
\end{array} \right).}
\ee
\end{prop}
{\bf Proof.}
Dropping the symmetrisation in (\ref{killing}) implies the existence of $\mu$
such that $\nabla_a K_b=\mu\epsilon_{ab}$. Differentiating this equation
covariantly, skew-symmetrising over all indices
and using  the curvature decomposition (\ref{formula_for_rho})
together with the Bianchi identity yields the statement of the Proposition.
\koniec
The connection $D$ is related to the standard tractor connection in projective differential geometry (see e.g. \cite{EM}). In the proof of Theorem
\ref{theo_1} we shall find the integrability conditions for the existence of parallel sections of this connection. This will lead to a set of invariants
of an affine connection $\nabla$.

{\bf Proof of Theorem \ref{theo_1}.}
The integrability conditions 
$(\nabla_a\nabla_b-\nabla_b\nabla_a)\mu=0$ give the algebraic
condition
\be
\label{obs}
F^aK_a+\beta\mu=0, \quad\mbox{where}\quad
F^a=\frac{1}{3}\epsilon^{ab}(L_b-\epsilon^{cd}\nabla_b B_{cd})
\ee
and $L_b\equiv \epsilon^{cd}\nabla_c\Rho_{db}$ is the Cotton tensor of $\nabla$. Geometrically the condition (\ref{obs}) means that the 
curvature of $D$ has rank at most one, and annihilates a parallel section of $D$. Applying $\nabla_a$ to
the condition (\ref{obs}), and using the vanishing of
(\ref{tractor_con}) leads to two more algebraic conditions
\be
{M_a}^{b}K_b+N_a\mu=0,
\label{obs2}
\ee
where
\[
{M_a}^b=\nabla_a F^b+\Big({\Rho^b}_a+\frac{1}{2}
{\delta^b}_a\beta\Big)\beta, \quad 
N_a=-F_a+\nabla_a\beta-\beta\theta_a.
\]
Multiplying the equation (\ref{obs}) by $2\theta_a$, and adding the resulting expression 
to (\ref{obs2}) results in
${M_a}^{b}\rightarrow {M_a}^{b}+2\theta_a F^b$ and
$N_a\rightarrow N_a+2\theta_a\beta$.
We can use this freedom to get rid of $\theta^a$ from the expressions
for $M$ and $N$. This yields
\be
\label{MandN}
{M_a}^b=\frac{1}{3}\epsilon^{bc}\epsilon^{de}
(\nabla_{a}Y_{dec}-\nabla_a\nabla_cB_{de})+\beta\;{\Rho^b}_a+\frac{1}{2}\beta^2{\delta^b}_a,
\quad N_a=-F_a+\epsilon^{bc}\nabla_a B_{bc},
\ee
where $Y_{cdb}=\nabla_{[c}\Rho_{d]b}$.
Therefore a parallel section $\Psi\equiv(K_1, K_2, \mu)^T$ 
of $D$ must satisfy
a system of three linear algebraic equations which we write
in a matrix form as
\be
\label{matrixm}
{\mathcal M}\Psi\equiv
\left(\begin{array}{ccc}
F^1 & F^2 & \beta\\
{M_1}^{1} & {M_1}^{2} & N_1\\
{M_2}^{1} & {M_2}^{2} & N_2  
\end{array} \right)
\left(\begin{array}{c} K_1 \\ 
K_2\\
\mu
\end{array} \right)=0.
\ee
A necessary condition for the existence of a non--zero parallel section $\Psi$ is therefore the vanishing of the determinant of 
the matrix ${\mathcal M}$. This gives the first obstruction which we write
as a vanishing of the relative scalar invariant
\be
\label{i_n}
I_N=\epsilon_{cd}\epsilon^{be}{M_e}^c\Big(N_bF^d-\frac{1}{2}\beta {M_b}^d\Big).
\ee
This invariant has weight $-5$: if we replace $\epsilon_{ab}$ by
$e^{f}\epsilon_{ab}$, where $f:\Sigma\rightarrow \R$, then
$I_N\rightarrow e^{-5f} I_N$. Thus $I_N\otimes(\epsilon_{ab}dX^a\wedge dX^b)^{\otimes 5}$ is an invariant.
The vanishing of $I_N$ is not sufficient for the existence of a non--zero parallel section. To assure sufficiency assume 
that $I_N=0$. Rewrite (\ref{obs}) and (\ref{obs2}) as
\[
V^{\alpha}\Psi_\alpha=0,\quad (D_a V^{\alpha})\Psi_\alpha=0,\quad \alpha=1, \dots, 3
\]
where $V=(F^1, F^2, \beta)$ in the formula above is a section of the dual bundle $E^*$,
and $D_a$ is the dual connection inherited from (\ref{tractor_con}). We continue differentiating, and 
adding the linear equations on $\Psi$. The Frobenius theorem tells
us that the process terminates once a differentiation
does not add any additional independent equations, as then the rank of the  
matrix of equations on $\Psi$ stabilises and does not grow. The space of parallel sections of $D$ has dimension equal to $3$ (the rank of the bundle $E$) minus the number of independent equations on $\Psi$. 
Therefore the sufficient condition for the existence of a Killing form assuming that $I_N=0$ is
\be
\label{6th}
{\mbox{rank}}{
\left(
\begin{array}{c}
{V}\\
{D_1 V}\\
{D_2 V}\\
{D_1D_1 V}\\
D_{(1}D_{2)}V\\
D_2D_2 V
\end{array}
\right )<3.}
\ee
If $I_N=0$ and  $V\neq 0$, then 
\[
cV+c_1D_1V+c_2D_2V=0,
\]
where $(c, c_1, c_2)$ are some functions on $U$, and $(c_1, c_2)$ are  not both zero. This implies that the term $D_{(1}D_{2)}V$ in 
(\ref{6th}) is a linear combination of all other terms, and can be disregarded. Now, suppose
that $D_1V=0$. Then (\ref{6th}) 
becomes $\mbox{det}(V, D_2 V, D_2 D_2 V)=0$. Equivalently, if $D_2V=0$ then
(\ref{6th})  becomes $\mbox{det}(V, D_1 V, D_1 D_1 V)=0$.
We conclude that (\ref{6th}) is equivalent to
\be
\label{i_s}
I_s=W_{abc}\equiv\mbox{det}(V, D_a V, D_{(b} D_{c)} V)=0.
\ee
as it is easy to show that the condition above implies (\ref{6th}). In fact vanishing of 
$W_{abc}$ gives just one independent condition:
If $c_2\neq 0$, then the sufficient condition is $W_{111}=0$, and if $c_1\neq 0$, then
it is $W_{222}=0$. The explicit tensor form of the obstruction $W$ is
\be
\label{w_tensor}
W_{acd}= F_b {M_a}^b V_{(cd)}-F_b U^b_{(cd)}N_a+\beta M_{ab}U^b_{(cd)}, \quad\mbox{where}
\ee
\begin{eqnarray*}
U^b_{ca}&=&\epsilon^{bd}\epsilon^{ef}[\frac{1}{3}(\nabla_c\nabla_a
Y_{efd}-\nabla_c\nabla_a\nabla_d B_{ef})+\nabla_c( B_{ef}P_{da})]\\
&+&\frac{1}{2}\epsilon^{ef}\epsilon^{gh}\nabla_c(B_{ef}
B_{gh})\delta^b_a+\epsilon^{bd}N_a(P_{dc}+\frac{1}{2}\beta \epsilon_{cd}),\quad\mbox{and}\\
V_{ca}&=&-M_{ac}-\frac{1}{3}\epsilon^{de}(\nabla_c\nabla_d
P_{ea}-\nabla_c\nabla_a B_{de}).
\end{eqnarray*}
\vskip2pt
We shall now consider the case when there exist precisely two independent solutions to the Killing equation
(\ref{killing}) (note that this situation does not arise if  
$\nabla$ is a Levi--Civita connection of some metric, 
as then the number of Killing vectors can be
$0, 1$ or $3$ - the last case being projectively flat). Therefore the rank of the matrix ${\mathcal{M}}$ in (\ref{matrixm}) is equal to one. We find that this can happens  if and only if  $\beta\neq 0$
and
\be
\label{tensor_obs}
{T_a}^{b}=0, \quad\mbox{where}\quad
{T_a}^{b}\equiv N_aF^b-\beta {M_a}^{b}.
\ee
This condition
guarantees the vanishing of all two-by-two minors of ${\mathcal M}$.
\vskip2pt
Finally, there exist three independent parallel sections of $D$ iff
the curvature of $D$ vanishes, or equivalently if the matrix 
${\mathcal M}$ vanishes. This condition is equivalent to the projective
flatness of the connection $\nabla$ together with the condition $\beta=0$.
\koniec
{\bf Remarks}
\begin{itemize}
\item
If the connection $\nabla$ is special (i.e. the Ricci tensor is symmetric, or equivalently $\beta=0$) then $I_N=-3^{-3}\nu_5$, where 
\[
\nu_5\equiv L^aL^b\nabla_aL_b
\]
is the Liouville projective invariant \cite{liouville, BDE}, and the indices are rised with a parallel
volume form.  Note that, unlike $\nu_5$, the obstruction
$I_N$ is not invariant under the projective changes of connection (see 
eq. (\ref{projective_change}) in \S\ref{sec_hydro}). The sufficient condition (\ref{i_s}) is then equivalent 
to the vanishing of the invariant $w_1$ constructed by Liouville for second order ODEs in 
\cite{liouville}.
\item
Theorem \ref{theo_1} generalises a well known characterisation of metrics 
which admit a Killing vector as those with functionally dependent scalar 
invariants. See \cite{Kruglikov} or \cite{MD} where a 3 by 3 matrix
analogous to ${\mathcal M}$ has been constructed. In this case
$N=-F= \frac{1}{3}*dR$, where $R$ is the scalar curvature, and $*$ is the Hodge operator of the metric $g$. The invariant (\ref{i_n}) reduces 
to\footnote{The prolongation procedure in \cite{MD} has been carried over in the Riemannian case. The additional subtlety in the 
Lorentzian signature arises if  $\nabla R$ is a non--zero null vector. 
We claim that no non-zero Killing vectors exist in this case. 
To see it, assume that a Lorentzian metric admits a Killing vector $K$. If $K$ is null, then the metric is flat with $R=0$.
Otherwise it  can locally be put in 
the form $dY^2-f(Y)^2dX^2$ for some $f=f(Y)$. Imposing the condition 
$|\nabla R|^2\equiv 0$ leads to $R=$ const.}
\[
I_N:=*\frac{1}{432}dR\wedge d(|\nabla R|^2).
\]
\item Any affine connection $\nabla$ on $\Sigma$ corresponds to a family of 
neutral signature  
anti--self--dual Einstein 
metrics on a certain rank-$2$ affine bundle $M\rightarrow \Sigma$, given by 
\cite{DM}
\[
G=\left(d\xi_a-\left(\Gamma_{ab}^c \xi_c-\Lambda \xi_a\xi_b- 
\Lambda^{-1}\Rho_{ba}\right)d X^b\right)\odot d X^a,
\]
where $\xi_a$ are local coordinates on the fibres of $M$, and
$-24\Lambda$ is the Ricci scalar. These metrics admit
a linear first integral iff $\nabla$ admits a projective vector 
field.
\end{itemize}
In the metric case, a Levi--Civita connection can not admit precisely two
local
linear first integrals, as $\beta$ (which is proportional to 
the skew part of the Ricci tensor) 
vanishes.
In the following Proposition we shall explicitly 
find all local normal forms of non--metric affine connections which admit two first integrals.
\begin{prop}
\label{prop_rank1}
Let $\nabla$ be an affine connection on a surface $\Sigma$ which admits exactly
two non--proportional linear first integrals which are independent at some point
$p\in \Sigma$. Local 
coordinates $(X, Y)$ can be chosen on an open set $U\subset \Sigma$ containing $p$ 
such that 
\be
\label{class_2}
\Gamma_{12}^1=\Gamma_{21}^1= \frac{c}{2}, \quad
\Gamma_{11}^2=\frac{P_X}{Q}, \quad \Gamma_{12}^2=  \Gamma_{21}^2 =\frac{P_Y+Q_X-cP}{2Q}, \quad
\Gamma_{22}^2=\frac{Q_Y}{Q},
\ee
and all other components vanish, where $c$ is a constant equal to $0$ or $1$,
and $(P, Q)$ are arbitrary functions of $(X, Y)$.
\end{prop}
{\bf Proof.}
Let the one--forms $K$ and $L$ be two solutions to the Killing equation.
If $K$ is closed, then there exist local coordinates $(X, Y)$ on $U$ such
that $K=dX$, and the corresponding first integral is $\dot{X}$. Therefore
$\ddot{X}=0$ and the connection components $\Gamma^1_{ab}$ vanish. Let the 
second solution of the Killing equation be of the form $L=PdX+QdY$ for some functions
$(P, Q)$. Imposing
\[
\frac{d}{d\tau}(P\dot{X}+Q\dot{Y})=0
\]
yields the non--zero components of the connection
given by (\ref{class_2}) with $c=0$.
If $dK\neq 0$, then coordinates $(X, Y)$ can be chosen so that $K=e^Y dX$.
The condition $d/d\tau(e^Y\dot{X})=0$ gives 
 $\Gamma^{1}_{12}=1/2$. Imposing the existence of the second integral
$(P\dot{X}+Q\dot{Y})$ yields the connection (\ref{class_2}) with $c=1$.
\koniec
Note that in both cases the ODEs for the unparametrised geodesics 
also admit a first integral,
given by $e^{-cY}({P+Y'Q})$, where $'=d/dX$. Conversely
if a 2nd order ODE cubic in $Y'$ representing
projective equivalence class $[\nabla]$
of affine connections admits a first integral linear in $Y'$, then $[\nabla]$ contains a connection
of the form (\ref{class_2}) with $c=0$. To see it consider a second order  ODE
of the form $(P+Y'Q)'=0$, where $(P, Q)$ are arbitrary functions of $(X, Y)$ 
and write it in the form
\be
\label{projective_ode}
Y''=\Gamma_{22}^1 (Y')^3+(2\Gamma_{12}^1-\Gamma_{22}^2) (Y')^2+(\Gamma_{11}^1-
2\Gamma_{12}^2)Y'-\Gamma_{11}^2.
\ee
Equation (\ref{projective_ode}) arises from eliminating the affine parameter 
$\tau$ between the two ODEs (\ref{flow}). Thus its integral curves are unparametrised geodesics of the affine connection $\nabla$.

\section{Hamiltonian systems of hydrodynamic type}
\label{sec_hydro}
An $n$--component $(1+1)$ system of hydrodynamic type has the form
$\p_t u^a={v^a}_b(u)\p_xu^b$, where $u^a=u^a(x, t)$ and $a, b=1, \dots, n$. From now on 
we shall assume that $n=2$ and that the matrix $v$ is diagonalisable at
some point with distinct eigenvalues, in which case there always exists
(in a neighbourhood of this point)
two distinct functions (called the Riemann invariants) $X^1$ and $X^2$ of $(u^1, u^2)$
such that the system is diagonal, i.e. takes the form 
(\ref{hydro_system})
for some $\lambda^a(X^b)$ and can be linearised by a hodograph transformation 
interchanging dependent $(X^1, X^2)$ and independent $(x, t)$ coordinates.

The hydrodynamic type system is said to admit a local 
Hamiltonian formulation with a
Hamiltonian of hydrodynamic type \cite{dubrovin, Fer}, if there exists a functional
$H[X^1, X^2]=\int_\R {\mathcal H}(X^1, X^2) dx$, where the density 
${\mathcal H}$ does not depend on
the derivatives of $X^a$ and  such that  (\ref{hami}) holds
for some functions $g^{ab}(X)$ and $b^{ab}_c(X)$.
If the matrix $g^{ab}$  is non--degenerate, then the Poisson bracket
\[
\{F, G\}=\int_\R \frac{\delta F}{\delta X^a} 
\Big(g^{ab}\frac{\p}{\p x}+{ b^{ab}_c}\frac{\p X^c}{\p x}\Big)
\frac{\delta G}{\delta X^b}dx 
\]
is skew--symmetric if $g^{ab}$ is symmetric and the metric
$g=g_{ab} dX^a dX^b$, where $g_{ab}g^{bc}={\delta_a}^c$ is parallel with respect to the connection with Christoffel symbols $\gamma_{ab}^c$ defined by $b^{ab}_c=-g^{ad}{\gamma^b}_{dc}$. The Jacobi identity then holds iff
the metric $g$ is flat, and the connection  defined by 
$\gamma_{ab}^c$ is torsion--free. The hydrodynamic type systems which 
admit a Hamiltonian of hydrodynamic type possess infinitely many
Poisson commuting first integrals, and are integrable in the Arnold--Liouville 
sense \cite{tsarev}.

{\bf Proof of Theorem \ref{theo_2}.}
It was shown in \cite{dubrovin} that a hydrodynamic 
type system in Riemann invariants is Hamiltonian in the sense defined above
if and only if there exists a flat diagonal metric
\be
\label{flat_met}
g=k^{-1} d(X^1)^2+f^{-1} d(X^2)^2
\ee
on a surface $U$ with local coordinates $(X^1, X^2)$
such that
\be
\label{eq1}
\p_2 k+2Ak=0, \quad\p_1 f +2B f=0,
\ee
where $f, k$ are functions of $(X^1, X^2)$, and $(A, B)$ are given by
(\ref{connection_{AB}}). Flatness of the metric $g$ yields
\be
\label{eq2}
(\p_2 A+A^2)f+(\p_1 B+B^2)k+\frac{1}{2}A\p_2 f+\frac{1}{2}B\p_1 k=0.
\ee
We verify that  equations (\ref{eq1}) and (\ref{eq2}) are equivalent
to  the Killing equations
(\ref{killing})
for an affine torsion--free connection $\nabla$ on $U$ defined by
(\ref{connection_{AB}}) where $K_1=Af, K_2=Bk$.
\koniec
Computing the relative invariants $I_N$ and $I_S$ gives explicit but 
complicated (albeit perfectly manageable by MAPLE) obstructions given in terms
of $(\lambda^1, \lambda^2)$ and their derivatives of order up to $6$.
These obstructions, together with the tensor (\ref{tensor_obs}) and the 
Cotton tensor of $\nabla$
characterise
Hamiltonian, bi-Hamiltonian and tri-Hamiltonian systems of hydrodynamic type.
The tri-Hamiltonian systems have been previously characterised by Ferapontov
in \cite{Fer} in terms of  two differential forms he called $\omega$ and 
$\Omega$. We shall now show how Ferapontov's formalism relates to our
connection (\ref{connection_{AB}}). We shall find that 
$\Omega$ is proportional to the skew-symmetric part of the Ricci tensor
of $\nabla$, and $\omega$ is the volume form of  the (generically) unique Lorentzian 
metric on $U$ which shares its unparametrised geodesics with $\nabla$.

\vskip5pt

We say that a symmetric affine connection $\nabla$ is {\em metric}, iff it is 
the Levi--Civita connection of some (pseudo)--Riemannian metric. An affine connection
$\nabla$ 
is {\em metrisable} iff it shares its unparametrised geodesic with some 
metric connection. Thus in the metrisable case there exists a one--form
$\Upsilon$ and a metric $h$ such that the Levi--Civita connection of $h$ is 
given by
\be
\label{projective_change}
{{\Gamma}^a}_{bc}+\delta^a_b\Upsilon_c+
\delta^a_c\Upsilon_b,
\ee
where ${{\Gamma}^a}_{bc}$ are the Christoffel  symbols of $\nabla$.
Not all affine connections on a surface are metrisable. The necessary and sufficient conditions for metrisability have been found in \cite{BDE}. 
\begin{prop}
The connection (\ref{connection_{AB}})  from Theorem \ref{theo_2}
is generically not metric but is metrisable by the 
metric
\be
\label{metric_h}
h=AB (dX^1)\odot (dX^2).
\ee
\end{prop}
{\bf Proof.} 
The connection is generically not metric, as its Ricci tensor $R_{ab}$ is in general 
not symmetric. The skew part of $R_{ab}$  is given by
\be
\label{upsilon}
(R_{21}-R_{12})dX^1\wedge dX^2=3d\Upsilon,\quad\mbox{where}\quad
\Upsilon=\Big(\frac{1}{2}\p_1 \ln{B}+B\Big)dX^1+
\Big(\frac{1}{2}\p_2 \ln{A}+A\Big)dX^2.
\ee
The unparametrised geodesics of this connection are integral curves of a 2nd order ODE
\be
\label{ode}
Y''=(\p_X Z)Y'-(\p_Y Z) (Y')^2, \quad \mbox{where}\quad Z=\ln{(AB)},
\ee
and $(X^1, X^2)=(X, Y)$. The ODE (\ref{ode}) is also the equation for 
unparametrised geodesics of the pseudo-Riemannian metric (\ref{metric_h})
(it can be found directly by solving the metricity equations as in \cite{CD16}).
The Levi--Civita connection of $h$ is given by (\ref{projective_change}),
where $\Upsilon$ is given by (\ref{upsilon}). Therefore $\nabla$ is projectively equivalent to a metric connection.
\koniec
{\bf Remarks}
\begin{itemize}
\item The pseudo--Riemannian metric (\ref{metric_h}) 
depends only on the product $AB$, so the transformation
$(A\rightarrow \gamma A, B \rightarrow \gamma^{-1}B)$,
where $\gamma=\gamma(X^a)$ is a non--vanishing function,
does not change 
unparametrised geodesics. It corresponds to a projective change of connection 
(\ref{projective_change}) by a one--form
\[
\Upsilon=
\Big((1-\gamma^{-1})B+\frac{1}{2}\p_1\ln{\gamma}\Big)dX^1
+\Big((1-\gamma)A-\frac{1}{2}\p_2\ln{\gamma}\Big)dX^2.
\]
This transformation can be used to set $R_{[ab]}$ to zero, but it does not 
preserve  (\ref{hydro_system}).
\item As the Ricci tensor $R_{ab}$ is in general not symmetric,
the connection (\ref{connection_{AB}}) does not admit a volume form on $\Sigma$ which is parallel w.r.t $\nabla$. Therefore
the Killing equations (\ref{killing}) do not imply the existence of a Killing vector for the metric $h$.
\item
The two--form $\Omega$ in Theorem $9$ of \cite{Fer} equals $2d\Upsilon$, while
$\omega$ in \cite{Fer} is given by the volume form of $h$.
In the tri--Hamiltonian case the connection $\nabla$ is projectively flat.
Equivalently  the metric (\ref{metric_h}) has constant Gaussian curvature, 
i.e. 
\be
\label{liouville}
(AB)^{-1}\p_1\p_2\ln{(AB)}=\mbox{const}.
\ee This is the Liouville
equation from Section 5 in \cite{Fer}.
\item If $n\geq 3$, there is always a discrepancy between the number of equation
for a Killing tensor of any given rank and a number of conditions for a HT system to admit a Hamiltonian formulation. Therefore Theorem \ref{theo_2} does not generalise to higher dimensions in any 
straightforward way.
\end{itemize}
\section{Examples}
\label{sec_ex}
In the examples below we set $X^1=X, X^2=Y$.
\subsection*{Example 1} Consider an affine connection (\ref{connection_{AB}}) corresponding to a system of hydrodynamic type with 
\[
A=cX+Y,  \quad B=X+cY,\quad \mbox{where}\quad c=\mbox{const.}
\]
This connection admits a parallel volume form iff $c=0$ or $c=1$.
If $c=0$ then the connection is projectively flat, and so the
system of hydrodynamic type is tri-hamiltonian.
Calculating the obstruction (\ref{tensor_obs}) yields
\[
T=
\frac{8c^2(c^2-9)}{9(cX+Y)^3(X+cY)^3}\Big(
dY\otimes\p_Y-
dX\otimes\p_X+
\frac{X+cY}{cX+Y}
dY\otimes \p_X-  \frac{cX+Y}{X+cY} dX\otimes \p_Y\Big).
\]
Therefore, if $c=3$ or $c=-3$ then the connection admits 
precisely two linear first integrals, so the system is bi-Hamiltonian. 
Finally for any $c$ not equal to $0, \pm 3$ the system admits a unique 
Hamiltonian.
\subsection*{Example 2} 
One dimensional non--linear elastic medium is governed by the system of 
PDEs \cite{olver,sheftel}
\[
u_t=h^2(v)v_x, \quad v_t=u_x,
\]
where $h(v)$ is a function characterising the type of fluid. This system is Hamiltonian with ${\mathcal H}=u^2/2+F(v)$, where $F''=h^2$. We find
the Riemann invariants $(X, Y)$ such that
\[
u=X+Y, \quad  v=G(X-Y),\quad\mbox{where} \quad G'h(G)=1, \quad\mbox{and}\quad
\lambda^1=-\lambda^2=\frac{1}{G'}.
\]
Therefore
$
A=-B=-G''/(2G')
$
and we find $\beta=0$ and so the Ricci tensor of the associated connection (\ref{connection_{AB}})
is symmetric. In particular, Theorem \ref{theo_1} implies that the system can not admit precisely two Hamiltonian structures.

 The projective flatness (\ref{liouville}) of the connection 
(\ref{connection_{AB}})
reduces to $(\ln{A^2})''=\mbox{const}.A^2$ which can be solved explicitly, 
and leads to a four-parameter family of tri-Hamiltonian systems. The singular solution $A=1/(2z)$ corresponds to the Toda equation
$v_{tt}=(\ln{v})_{xx}$.
\subsection*{Example 3} We consider
the system of hydrodynamic type (\ref{hydro_system}) with
\[
\lambda^1= -\lambda^2=(X-Y)^n (X+Y)^m.
\]
Examining the conditions of Theorem \ref{theo_1} for the resulting 
connection (\ref{connection_{AB}})
we find that this system is always bi-Hamiltonian. It is 
tri-Hamiltonian iff $nm(n^2-m^2)=0$.
\subsection*{Example 4. Frobenius manifolds}
In this example we shall consider Hamiltonian systems of hydrodynamic type which arise
from two--dimensional Frobenius manifolds.
Recall \cite{dubrovin_frob,hitchin_frob} that $\Sigma$ is a two--dimensional 
Frobenius manifold if the tangent space to $\Sigma$ at each point admits a structure of a 
 commutative algebra ${\bf A}$ with a unity ${\bf e}$, and symmetric 
tensors fields
$C\in C^{\infty}(\mbox{Sym}^3(T^*\Sigma))$  and 
$\eta \in C^{\infty}(\mbox{Sym}^2(T^*\Sigma))$ such that locally, 
in $U\subset \Sigma$, 
there exist a coordinate system $u^a=(u, v)$ and a function $F:U\rightarrow \R$ 
where
\[
C=\frac{\p^3 F}{\p u^a \p u^b \p u^c} du^a du^b du^c, \quad {\bf e} =\frac{\p}{\p u^1},\quad
\eta=\frac{\p^3 F}{\p u^1 \p u^a \p u^b}du^a du^b. 
\]
The non--degenerate symmetric form $\eta$ is
a flat
(pseudo) Riemannian metric with constant coefficents 
on $\Sigma$ such that ${\bf e}$ is covariantly constant, 
and ${C^a}_{bc}:=\eta^{ad}C_{bcd}$
are the structure constants  for ${\bf A}$. Moreover there exists an Euler vector field
${\bf E}$ such that ${\mathcal L}_{\bf E} {\bf e}=-{\bf e}$, and
${\mathcal L}_{\bf E} C=(m+3)C, 
{\mathcal L}_{\bf E} \eta=(m+2)\eta$ for some constant $m$. 

In dimensions higher that two the function
$F$ must satisfy a non--linear PDE resulting from the associativity conditions.
In two dimensions the associativity always holds, and $F$ can be found
only from the homogeneity condition. If we assume that the identity vector field
${\bf e}$ is null with respect to the metric $\eta$, then we can set
$\eta=du\odot dv$, and find \cite{dubrovin_frob,hitchin_frob} that
$F(u, v)=\frac{1}{2}u^2v+f(v)$, where $f$ (which we assume not to be zero)
is given by one of the four expressions
\[
f=v^k,\; k\neq 0, 2, \quad f=v^2\ln v, \quad f=\ln v, \quad f=e^{2v}.
\]
In all cases the corresponding Hamiltonian system of hydrodynamic type  is
\be
\label{f_system}
u_t=f'''(v) v_x, \quad v_t=u_x.
\ee
The characteristic velocities are $\lambda^1=-\lambda^2=\lambda\equiv\sqrt{f'''(v)}$,
and the Riemann invariants $X^a=(X, Y)$ are
\[
X=u+\int \sqrt{f'''(v)} dv, \quad Y=u-\int \sqrt{f'''(v)} dv.
\]
The corresponding affine connection (\ref{connection_{AB}}) is projectivelly flat 
(the Cotton tensor $\nabla_{[a}\Rho_{b]c}$ vanishes), and special 
(the Ricci tensor is symmetric). Therefore 
Theorem \ref{theo_1} implies that the system (\ref{f_system}) is tri--Hamiltonian.
The corresponding three--parameter family of flat metrics (\ref{flat_met}) is 
\be
g(c_1, c_2, c_3)=
\lambda^{-1}\Big(\frac{dX^2}{c_1+c_2X+c_3X^2}-\frac{dY^2}{c_1+c_2Y+c_3Y^2}\Big),\quad
\lambda\equiv\sqrt{f'''(v)},
\ee
where $(c_1, c_2, c_3)$ are arbitrary constants not all zero. We find
that $\eta\equiv g(1, 0, 0)$ is the flat metric in the definition of the Frobenius manifold. The second metric $I\equiv g(0, 1, 0)$ is the so-called intersection
form (see \cite{dubrovin_frob}). The third metric
is  $J\equiv g(0, 0, 1)$. It can be constructed directly from $\eta$ and $I$ as
$J_{ab}=I_{ac}I_{bd}\eta^{cd}$ in agreement with \cite{pavlov,romano}. 
It would be interesting to analyse the existence of Hamiltonians for  HT type systems
arising from submanifolds of Frobenius manifolds \cite{Strachan}.
\subsection*{Example 5. Zoll connections}
Recall that a Riemannian metric $h$ on a surface $\Sigma$ is {\em Zoll} if all
geodesics are simple closed curves of equal length. A two--dimensional sphere
admits a family of axisymmetric Zoll metrics given by
\be
\label{zoll_1}
h=(F(X)-1)^2dX^2+\sin^2{X}dY^2,
\ee
where $(X, Y)$ are spherical polar coorinates on $\Sigma=S^2$, and
$F:[0, \pi]\rightarrow [0, 1]$ is any function such that $F(0)=F(\pi)=0$ and
$F(\pi-X)=-F(X)$. 
A projective structure $[\nabla]$ on $\Sigma$ is {\em Zoll} if its unparametrised geodesics are  simple closed curves. The general projective structure admitting a projective vector field, and close to the flat structure of the round sphere 
is given by the second order ODE \cite{zoll}
\be
\label{zoll_2}
Y''=A_3 (Y')^3+A_2(Y')^2+A_1Y',\quad\mbox{where}
\ee
\[
A_1=\frac{F'}{F-1}-2\cot{X}, \quad A_2=\frac{H'\sin{X}\cos{X} -2H}{\cos{X}(F-1)},
\quad A_3=-\frac{(H^2+1)\sin{X}\cos{X}}{(F-1)^2},
\]
where $F=F(X)$ is as before, and $H=H(X)$ satisfies $H(0)=H(\pi)=H(\pi/2)=0$,
and $H(\pi-X)=H(X)$. The metric case (\ref{zoll_1}) arises 
if $H=0$. A general 
connection $\nabla$ in this projective class with $\beta\neq 0$ will not admit
even a single first integral. We use Theorem \ref{theo_1} together with
(\ref{projective_ode}) to verify that the following choice of the 
representative connection 
\be
\label{zoll_3}
\Gamma_{11}^1=A_1, \quad \Gamma_{22}^1=A_3, \quad
\Gamma_{12}^1=\Gamma_{21}^1=\frac{1}
{2}A_2 
\ee
admits a first integral for any $F$ and $H$. 
To find a (necessarily non-metric) Zoll connection 
with precisely two linear first integrals we 
use Proposition (\ref{prop_rank1}) and match the connection 
(\ref{zoll_3}) with the connection
(\ref{class_2}) (with the roles of $X$ and $Y$ reversed). 
This, for any given $H$, leads to a one--parameter 
family of examples
\[
F=1+c(H^2+1)\cot{X}
\]
which does not satisfy the boundary conditions.
The existence of a non-metric Zoll structure on $S^2$ with precisely two first 
integrals is an interesting open problem.

\end{document}